\documentclass[francais]{article}

\usepackage{graphicx}
\usepackage{amsmath}
\usepackage{amsfonts}
\usepackage{amssymb}

\newtheorem{theorem}{Th\'eor\`eme}

\newtheorem{proposition}{Proposition}
\newtheorem{corollary}{Corollaire}
\newtheorem{definition}{D\'efinition\rm}
\newtheorem{remark}{Remarque}

\newenvironment{proof}[1][D\'emonstration]{\textbf{#1.} }{\ \rule{0.5em}{0.5em}}

\begin{document}

\title{Sur les alg\`{e}bres de Lie quasi-filiformes admettant un tore de d\'{e}rivations
}
\author{L. Garc\'{\i}a Vergnolle\footnote{e-mail:
lucigarcia@mat.ucm.es}\\
Dpto. Geometr\'{\i}a y Topolog\'{\i}a,\\
              Facultad de Ciencias Matem\'{a}ticas U.C.M.\\
              Plaza de Ciencias 3, 28040 Madrid
              Espagne\\
}

\date{ }

\maketitle

\begin{abstract}
Dans ce travail on d\'{e}crit les alg\`{e}bres de Lie quasi-filiformes de rang non nul. De plus, on rappelle et corrige la classification des alg\`{e}bres de Lie filiformes admettant un tore de d\'{e}rivations, ainsi que la liste des alg\`{e}bres gradu\'{e}es naturellement et quasi-filiformes.\\

\end{abstract}

\section{Introduction}

Soit $\frak{g}$ une alg\`{e}bre de Lie nilpotente de dimension $n$ et nilindice $m$, elle est naturellement filtr\'{e}e par la suite centrale descendante:\\
\begin{align*}
\frak{g}_{1}  &  =\frak{g}\supseteq \frak{g}_{2}=\left[
\frak{g},\frak{g}\right]\supseteq \frak{g}_{3}=\left[
\frak{g}_{2},\frak{g}\right]  \supseteq...\supseteq
\frak{g}_{k+1}=\left[  \frak{g}_{k},\frak{g}\right]  \supseteq...
\supseteq \frak{g}_{m+1}=\left\{ 0 \right\}
\end{align*}
On peut alors associer une alg\`{e}bre de Lie gradu\'{e}e, not\'{e}e par $gr(\frak{g})$, et d\'{e}finie par:
$$
{\rm gr}\frak{g} \; = \; \sum_{i=1}^{m} \, \frac{\frak{g}_{i}}{\frak{g}_{i+1}} \; = \; \sum_{i=1}^{m} W_{i}
$$
o\`{u} le crochet est donn\'{e} par:
$$[X+\frak{g}_{i+1},Y+\frak{g}_{j+1}]=[X,Y]+\frak{g}_{i+j+1}\quad \forall X \in \frak{g}_{i}\quad \forall Y \in \frak{g}_{j}$$
On dit que $\frak{g}$ est une alg\`{e}bre du type $\{p_{1},\dots,p_{m}\}$ si $\dim \frac{\frak{g}_{i}}{\frak{g}_{i+1}}=p_{i}$. Remarquons que $gr(\frak{g})$ est du m\^{e}me type que $\frak{g}$.\\
\begin{definition}
Une alg\`{e}bre de Lie nilpotente de dimension $n$ est dite filiforme si son nilindice est $n-1$, c'est-\`{a}-dire si elle est du type $\{2,1,1,\dots,1\}$
\end{definition}
On en d\'{e}duit que l'alg\`{e}bre gradu\'{e}e d'une alg\`{e}bre filiforme est aussi filiforme.
\begin{definition}
Une alg\`{e}bre $\frak{g}$ est gradu\'{e}e naturellement quand elle est isomorphe \`{a} ${\rm gr}\,\frak{g}$.
\end{definition}
Une alg\`{e}bre de Lie $\frak{g}$ gradu\'{e}e naturellement peut alors se d\'{e}composer de la forme $\frak{g}=W_{1}\oplus W_{2}\oplus \dots \oplus W_{m}$ o\`{u} les  $W_{i}=\frac{\frak{g}_{i}}{\frak{g}_{i+1}}$ v\'{e}rifient $[W_{i},W_{j}]=W_{i+j}$ pour $i+j\le m$.
\par
La classification des alg\`{e}bres de Lie filiformes gradu\'{e}es naturellement est due \`{a} Vergne \cite{Vergne}.
\begin{theorem}
Toute alg\`{e}bre de Lie filiforme et gradu\'{e}e naturellement est isomorphe \`{a} une des alg\`{e}bres suivantes:
\begin{enumerate}
\item $L_{n}\:(n\geq3)$ d\'{e}finie dans la base $\{X_{0},X_{1},\dots,X_{n-1}\}$ par:
$$
\lbrack X_{0},X_{i}]=X_{i+1},\;1\leq i\leq n-2.
$$
(les crochets non \'{e}crits \'{e}tant nuls, except\'{e}s ceux d\'{e}coulant de l'antisym\'{e}trie)
\item $Q_{n}\:(n=2m,\,m\geq3)$ d\'{e}finie dans la base $\{X_{0},X_{1},\dots,X_{n-1}\}$ par:
\begin{align*}
\lbrack X_{0},X_{i}]  & =X_{i+1},\;1\leq i\leq n-3,\\
\lbrack X_{j},X_{n-j-1}]  & =(-1)^{j-1}X_{n-1},\;1\leq j\leq m-1.
\end{align*}
\end{enumerate}
En cons\'{e}quence, sauf isomorphisme,il existe deux alg\`{e}bres filiformes gradu\'{e}es naturellement de dimension paire y une seule alg\`{e}bre filiforme gradu\'{e}e naturellement de dimension impaire.
\end{theorem}

\begin{proposition}
\cite{Goze} Soit $\frak{g}$ une alg\`{e}bre de Lie filiforme de dimension $n$ admettant une d\'{e}rivation diagonale $f$ non nulle. Il existe alors une base de $\frak{g}$, $\{Y_{0},Y_{1},\dots,Y_{n-1}\}$, form\'{e}e de vecteurs propres de $f$ dont les crochets v\'{e}rifient l'un des cas suivants:
\begin{enumerate}
\item $\frak{g}=L_{n},\: n\geq3$
$$
\lbrack Y_{0},Y_{i}\rbrack=Y_{i+1}, \quad 1 \le i \le n-2
$$
\item $\frak{g}=A_{n}^{k}(\alpha_{1},\dots,\alpha_{t-1}),\: n\geq3,\: t=\frac{n-k+1}{2},\: 2 \le k \le n-3$
$$
\begin{array}{ll}
\lbrack Y_{0},Y_{i}\rbrack=Y_{i+1}, & 1 \le i \le n-2,\\
\lbrack Y_{i},Y_{i+1}\rbrack=\alpha_{i} Y_{2i+k}, & 1 \le i \le t-1,\\
\lbrack Y_{i},Y_{j}\rbrack=a_{i,j} Y_{i+j+k-1}, & 1 \le i <j,\quad i+j \le n-k\\
\end{array}
$$
\item $\frak{g}=Q_{n},\: n=2m,\: m\geq 3$
$$
\begin{array}{ll}
\lbrack Y_{0},Y_{i}\rbrack=Y_{i+1}, & 1 \le i \le n-3,\\
\lbrack Y_{j},Y_{n-j-1}\rbrack  =(-1)^{j-1}Y_{n-1}, & 1\leq j\leq m-1\\
\end{array}
$$
\item $\frak{g}=B_{n}^{k}(\alpha_{1},\dots,\alpha_{t-1}) ,\: n=2m,\: m\geq 3 ,\: t=\frac{n-k}{2} ,\: 2 \le k \le n-3$
$$
\begin{array}{ll}
\lbrack Y_{0},Y_{i}\rbrack=Y_{i+1}, & 1 \le i \le n-3,\\
\lbrack Y_{j},Y_{n-j-1}\rbrack=(-1)^{j-1}Y_{n-1}, & 1\leq j\leq m-1,\\
\lbrack Y_{i},Y_{i+1}\rbrack=\alpha_{i} Y_{2i+k}, & 1 \le i \le t-1,\\
\lbrack Y_{i},Y_{j}\rbrack=a_{i,j} Y_{i+j+k-1}, & 1 \le i <j,\: i+j \le n-k-1\\
\end{array}
$$
\end{enumerate}
o\`{u} $(\alpha_{1},\dots,\alpha_{t-1})$ sont des param\`{e}tres v\'{e}rifiant les relations polynomiales d\'{e}coulant des identit\'{e}s de Jacobi et les contantes $a_{i,j}$ verifient le syst\`{e}me:
$$
\begin{array}{l}
a_{i,i}=0,\\
a_{i,i+1}=\alpha_{i},\\
a_{i,j}=a_{i+1,j}+a_{i,j+1}.
\end{array}
$$
De plus, si $\lambda_{0}$ et $\lambda_{1}$ repr\'{e}sentent les valeurs propres de $f$ associ\'{e}es respectivement aux vecteurs propres $Y_{0}$ et $Y_{1}$, alors $f$ prend, dans chacun des cas, la forme diagonale suivante:
\begin{enumerate}
\item $f\sim {\rm diag}(\lambda_{0},\lambda_{1},\lambda_{0}+\lambda_{1},2\lambda_{0}+\lambda_{1},\dots,(n-2)\lambda_{0}+\lambda_{1})$
\item $f\sim {\rm diag}(\lambda_{0},k\lambda_{0},(k+1)\lambda_{0},(k+2)\lambda_{0},\dots,(n+k-2)\lambda_{0})$
\item $f\sim {\rm diag}(\lambda_{0},\lambda_{1},\lambda_{0}+\lambda_{1},2\lambda_{0}+\lambda_{1},\dots,(n-3)\lambda_{0}+\lambda_{1},(n-3)\lambda_{0}+2\lambda_{1})$
\item $f\sim {\rm diag}(\lambda_{0},k\lambda_{0},(k+1)\lambda_{0},(k+2)\lambda_{0},\dots,(n+k-3)\lambda_{0},(n+2k-3)\lambda_{0})$
\end{enumerate}
\end{proposition}
\begin{remark}
Dans la r\'{e}f\'{e}rence \cite{Goze}, on consid\`{e}re aussi les alg\`{e}bres $C_{n}(\alpha_{1},\dots,\alpha_{m-2})$ o\`{u} $n=2m$, d\'{e}finies dans la base $\{Y_{0},Y_{1},\dots,Y_{n-1}\}$  par:
$$
\begin{array}{ll}
\lbrack Y_{0},Y_{i}\rbrack=Y_{i+1}, & 1 \le i \le n-3,\\
\lbrack Y_{i},Y_{n-i-1}\rbrack=(-1)^{i}Y_{n-1}, & 1\leq i\leq m-1,\\
\lbrack Y_{i},Y_{n-(2k+1)-i}\rbrack=(-1)^{i+1}\alpha_{k}Y_{n-1}, & 1\le i\le m-k-1 ,\:  1\le k\le m-2.\\
\end{array}
$$
Il existe un changement de variable tel que ces alg\`{e}bres adoptent la m\^{e}me forme avec $\alpha_{1}=\dots=\alpha_{m-2}=0$. En effet, pour \'{e}liminer le param\`{e}tre $\alpha_{1}$, il suffit de faire le changement de variable suivant:
$$
\begin{array}{ll}
Z_{0}=Y_{0}, & \\
Z_{1}=Y_{1}+\frac{\alpha_{1}}{2}Y_{3}, & \\
Z_{i+1}=\lbrack Z_{0},Z_{i}\rbrack, & 1 \le i \le n-3,\\
Z_{n-1}=Y_{n-1}, & \\
\end{array}
$$
et, de plus, si $\alpha_{j+1}=\dots=\alpha_{m-2}=0$, le changement de variable $Y_{i}\rightarrow Y_{i}+\frac{\alpha_{j}}{2}Y_{i+2j}$ pour $i=1,\dots,n-2-2j$,  nous permet annuler le param\`{e}tre $\alpha_{j}$.\\
On en conclut que les alg\`{e}bres $C_{n}$ sont isomorphes \`{a} $Q_{n}$; d\'{e}sormais les alg\`{e}bres de la forme $C_{n}$ ne sont pas de rang $1$ mais de rang $2$.
\end{remark}
\bigskip
Le but de ce travail est de donner la classification des alg\`{e}bres quasi-filiformes de rang non nul.

\section{Alg\`{e}bres de Lie quasi-filiformes et bases adapt\'{e}es}
De fa\c{c}on analogue \`{a} Vergne, ont \'{e}t\'{e} obtenues les alg\`{e}bres de Lie gradu\'{e}es naturellement et quasi-filiformes que l'on d\'{e}finit ci-dessous \cite{Gomez}.
\begin{definition}
Soit $\frak{g}$ une alg\`{e}bre nilpotente de dimension $n$, on dit que $\frak{g}$ es quasi-filiforme si son nil\'{\i}ndice est $n-2$.
\end{definition}

Si une alg\`{e}bre $\frak{g}$ est quasi-filiforme, en suivant les notations pr\'{e}c\'{e}dentes, il existe deux possibilit\'{e}s:
\begin{enumerate}
\item Soit $\frak{g}$ est du type $t_{1}=\{p_{1}=3,p_{2}=1,p_{3}=1,\dots,p_{n-2}=1\}$.
\item Soit $\frak{g}$ est du type $t_{r}=\{p_{1}=2,p_{2}=1,\dots,p_{r-1}=1,p_{r}=2,p_{r+1}=1,\dots,p_{n-2}=1\}$ o\`{u} $r\in \{2,\dots,n-2\}$.
\end{enumerate}

\begin{proposition}
\label{propo}
\cite{Gomez} Soit $\frak{g}$ une alg\`{e}bre de Lie quasi-filiforme gradu\'{e}e naturellement de dimension $n$ et du type $t_{r}$ o\`{u} $r\in \{1,\dots,n-2\}$. Il existe alors une base homog\`{e}ne $\{X_{0},X_{1},X_{2},\dots,X_{n-1}\}$ de $\frak{g}$ avec  $X_{0}$ et $X_{1}$ dans $W_{1}$, $X_{i}\in W_{i}$ pour $i\in \{2,\dots,n-2\}$ et $X_{n-1}\in W_{r}$ dans laquelle $\frak{g}$ est une des alg\`{e}bres d\'{e}crites ci-dessous.
\begin{enumerate}
\item $\frak{g}$ est du type $t_{1}$
\begin{enumerate}
\item $L_{n-1}\oplus \mathbb{C} \quad (n\ge 4)$
$$
\lbrack X_{0},X_{i}\rbrack=X_{i+1}, \quad 1 \le i \le n-3.
$$
\item $Q_{n-1}\oplus \mathbb{C} \quad (n\ge 7,\: n\;{\rm impair})$
$$
\begin{array}{ll}
\lbrack X_{0},X_{i}\rbrack=X_{i+1}, & 1 \le i \le n-4,\\
\lbrack X_{i},X_{n-i-2}\rbrack=(-1)^{i-1}X_{n-2}, & 1\leq i\leq \frac{n-3}{2}.\\
\end{array}
$$
\end{enumerate}
\item $\frak{g}$ est du type $t_{r}$ o\`{u} $r\in \{2,\dots,n-2\}$
\begin{enumerate}
\item $\frak{L_{n,r}}; \quad n\ge5, \;r\: {\rm impair}, \;3\le r\le 2[\frac{n-1}{2} ]-1$
$$
\begin{array}{ll}
\lbrack X_{0},X_{i} \rbrack=X_{i+1}, &  i=1,\dots,n-3\\
\lbrack X_{i},X_{r-i} \rbrack=(-1)^{i-1}X_{n-1}, &  i=1,\dots,\frac{r-1}{2}\\
\end{array}
$$
\item $\frak{Q_{n,r}}; \quad n\ge7,\;n\: {\rm impair}, \;r\: {\rm impair}, \;3\le r\le n-4$
$$
\begin{array}{ll}
\lbrack X_{0},X_{i} \rbrack=X_{i+1}, &  i=1,\dots,n-4\\
\lbrack X_{i},X_{r-i} \rbrack=(-1)^{i-1}X_{n-1}, &  i=1,\dots,\frac{r-1}{2}\\
\lbrack X_{i},X_{n-2-i} \rbrack=(-1)^{i-1}X_{n-2}, &  i=1,\dots,\frac{n-3}{2}\\
\end{array}
$$
\item $\frak{T_{n,n-4}};\quad n\ge7,\:n\; {\rm impair}$
$$
\begin{array}{ll}
\lbrack X_{0},X_{i} \rbrack=X_{i+1}, &  i=1,\dots,n-5\\
\lbrack X_{0},X_{n-3} \rbrack=X_{n-2},\\
\lbrack X_{0},X_{n-1} \rbrack=X_{n-3},\\
\lbrack X_{i},X_{n-4-i} \rbrack=(-1)^{i-1}X_{n-1}, &  i=1,\dots,\frac{n-5}{2}\\
\lbrack X_{i},X_{n-3-i} \rbrack=(-1)^{i-1} \frac{n-3-2i}{2} X_{n-3}, &  i=1,\dots,\frac{n-5}{2}\\
\lbrack X_{i},X_{n-2-i} \rbrack=(-1)^{i} (i-1) \frac{n-3-i}{2} X_{n-2}, &  i=2,\dots,\frac{n-3}{2}\\
\end{array}
$$
\item $\frak{T_{n,n-3}};\quad n\ge6,\:n \;{\rm pair}$
$$
\begin{array}{ll}
\lbrack X_{0},X_{i} \rbrack=X_{i+1}, &  i=1,\dots,n-4\\
\lbrack X_{0},X_{n-1} \rbrack=X_{n-2},\\
\lbrack X_{i},X_{n-3-i} \rbrack=(-1)^{i-1}X_{n-1}, &  i=1,\dots,\frac{n-4}{2}\\
\lbrack X_{i},X_{n-2-i} \rbrack=(-1)^{i-1} \frac{n-2-2i}{2} X_{n-2}, &  i=1,\dots,\frac{n-4}{2}\\
\end{array}
$$
\item $\frak{E_{9,5}^{1}}$
$$
\begin{array}{ll}
\lbrack X_{0},X_{i}\rbrack =X_{i+1}, & i=1,2,3,4,5,6 \\
\lbrack X_{0},X_{8}\rbrack =X_{6}, & \lbrack X_{2},X_{8}\rbrack =-3X_{7},\\
\lbrack X_{1},X_{4}\rbrack =X_{8}, & \lbrack X_{1},X_{5}\rbrack =2X_{6},\\
\lbrack X_{1},X_{6}\rbrack =3X_{7}, & \lbrack X_{2},X_{3}\rbrack =-X_{8},\\
\lbrack X_{2},X_{4}\rbrack =-X_{6}, & \lbrack X_{2},X_{5}\rbrack =-X_{7}.\\
\end{array}
$$
\item $\frak{E_{9,5}^{2}}$
$$
\begin{array}{ll}
\lbrack X_{0},X_{i}\rbrack =X_{i+1}, &i=1,2,3,4,5,6\\
\lbrack X_{0},X_{8}\rbrack =X_{6}, & \lbrack X_{2},X_{8}\rbrack =-X_{7},\\
\lbrack X_{1},X_{4}\rbrack =X_{8}, & \lbrack X_{1},X_{5}\rbrack =2X_{6},\\
\lbrack X_{1},X_{6}\rbrack =X_{7}, & \lbrack X_{2},X_{3}\rbrack =-X_{8},\\
\lbrack X_{2},X_{4}\rbrack =-X_{6}, & \lbrack X_{2},X_{5}\rbrack =X_{7},\\
\lbrack X_{3},X_{4}\rbrack =-2X_{7}. &\\
\end{array}
$$
\item $\frak{E_{9,5}^{3}}$
$$
\begin{array}{ll}
\lbrack X_{0},X_{i}\rbrack =X_{i+1}, &i=1,2,3,4,5,6\\
\lbrack X_{0},X_{8} =X_{6}, & \lbrack X_{1},X_{4}\rbrack =X_{8},\\
\lbrack X_{1},X_{5}\rbrack =2X_{6}, & \lbrack X_{2},X_{3}\rbrack =-X_{8},\\
\lbrack X_{2},X_{4}\rbrack =-X_{6}, & \lbrack X_{2},X_{5}\rbrack =2X_{7},\\
\lbrack X_{3},X_{4}\rbrack =-3X_{7}. &\\
\end{array}\\
$$
\item $\frak{E_{7,3}}$
$$
\begin{array}{ll}
\lbrack X_{0},X_{i}\rbrack =X_{i+1}, &i=1,2,3,4\\
\lbrack X_{0},X_{6}\rbrack =X_{4}, & \lbrack X_{2},X_{6}\rbrack =-X_{5},\\
\lbrack X_{1},X_{2}\rbrack =X_{6}, & \lbrack X_{1},X_{3}\rbrack =X_{4},\\
\lbrack X_{1},X_{4}\rbrack =X_{5}. &\\
\end{array}\\
$$
\end{enumerate}
\end{enumerate}

\end{proposition}
On en d\'{e}duit qu'il n'y a pas d'alg\`{e}bres quasi-filiformes du type $t_{2}$.
\begin{remark}
L'alg\`{e}bre  $\frak{E_{9,5}^{3}}$ ne se touve pas dans la classification \cite{Gomez}. En effet, quand on cherche les alg\`{e}bres de dimension $9$ ayant une d\'{e}rivation diagonale de la forme ${\rm diag}(1,1,2,3,4,5,6,7,5)$, on obtient facilement trois alg\`{e}bres sp\'{e}ciales, isomorphes respectivement \`{a} $\frak{E_{9,5}^{1}}$, $\frak{E_{9,5}^{2}}$ et $\frak{E_{9,5}^{3}}$. De plus, cette alg\`{e}bre est d'une importance consid\'{e}rable dans le probl\`{e}me de rigidit\'{e} \cite{AG} \cite{GR}.
\end{remark}
\begin{corollary}
Soit $\frak{g}$ une alg\`{e}bre de Lie quasi-filiforme gradu\'{e}e naturellement de dimension $n$ . Alors, ${\rm gr}\frak{g}$ est isomorphe \`{a} l'une des alg\`{e}bres de la proposition ~\ref{propo}, $C_{i,j}^{k}$ d\'{e}notent les constantes de structure de cette alg\`{e}bre dans la base o\`{u} on l'a d\'{e}finie. Il existe alors une base $\{X_{0},X_{1},X_{2},\dots,X_{n-1}\}$ de $\frak{g}$ telle que:
$$
\begin{array}{ll}
\lbrack X_{0},X_{i}\rbrack-\sum_{k=0}^{n-1}C_{0,i}^{k}X_{k}\in \frak{g}_{i+1}, & 1 \le i \le n-3 ,\\
\lbrack X_{i},X_{j}\rbrack-\sum_{k=0}^{n-1}C_{i,j}^{k}X_{k}\in \frak{g}_{i+j}, & 1 \le i<j \le n-2-i,\\
\lbrack X_{i},X_{n-1}\rbrack-\sum_{k=0}^{n-1}C_{i,n-1}^{k}X_{k}\in \frak{g}_{i+r}, & 1 \le i \le n-2-r,\\
\lbrack X_{i},X_{n-2}\rbrack=0, & 0 \le i \le n-1,\\
\end{array}
$$
et, de plus:
\begin{enumerate}
\item Si ${\rm gr}\frak{g}\simeq \frak{L_{n-1}}\oplus \mathbb{C}$,
$$
\lbrack X_{0},X_{i}\rbrack=X_{i+1}, \quad 1 \le i \le n-3,
$$
\item Si ${\rm gr}\frak{g}\simeq \frak{Q_{n-1}}\oplus \mathbb{C}$,
$$
\begin{array}{ll}
\lbrack X_{0},X_{i}\rbrack=X_{i+1}, & 1 \le i \le n-4,\\
\lbrack X_{1},X_{n-2}\rbrack=X_{n-2}, &\\
\end{array}
$$
\item Si ${\rm gr}\frak{g}\simeq \frak{L_{n,r}}$,
$$
\begin{array}{ll}
\lbrack X_{0},X_{i} \rbrack=X_{i+1}, &  i=1,\dots,n-3\\
\lbrack X_{1},X_{r-1} \rbrack=X_{n-1}, &\\
\end{array}
$$
\item Si ${\rm gr}\frak{g}\simeq \frak{Q_{n,r}}$,
$$
\begin{array}{ll}
\lbrack X_{0},X_{i} \rbrack=X_{i+1}, &  i=1,\dots,n-4\\
\lbrack X_{1},X_{r-1} \rbrack=X_{n-1}, &\\
\lbrack X_{1},X_{n-3} \rbrack=X_{n-2}, &\\
\end{array}
$$
\item Si ${\rm gr}\frak{g}\simeq \frak{T_{n,n-4}}$,
$$
\begin{array}{ll}
\lbrack X_{0},X_{i} \rbrack=X_{i+1}, &  i=1,\dots,n-5\\
\lbrack X_{0},X_{n-3} \rbrack=X_{n-2}, &\\
\lbrack X_{0},X_{n-1} \rbrack=X_{n-3}, &\\
\lbrack X_{1},X_{n-5} \rbrack=X_{n-1} &\\
\end{array}
$$
\item Si ${\rm gr}\frak{g}\simeq \frak{T_{n,n-3}}$,
$$
\begin{array}{ll}
\lbrack X_{0},X_{i} \rbrack=X_{i+1}, &  i=1,\dots,n-4\\
\lbrack X_{0},X_{n-1} \rbrack=X_{n-2}, &\\
\lbrack X_{1},X_{n-4} \rbrack=X_{n-1}, &\\
\end{array}
$$
\item Si ${\rm gr}\frak{g}\simeq \frak{E_{9,5}^{j}}$ avec $j\in\{1,2,3\}$,
$$
\begin{array}{ll}
\lbrack X_{0},X_{i}\rbrack =X_{i+1}, &i=1,2,3,4,5,6\\
\lbrack X_{1},X_{4}\rbrack =X_{8}, &\\
\end{array}\\
$$
\item Si ${\rm gr}\frak{g}\simeq \frak{E_{7,3}}$,
$$
\begin{array}{ll}
\lbrack X_{0},X_{i}\rbrack =X_{i+1}, &i=1,2,3,4\\
\lbrack X_{1},X_{2}\rbrack =X_{6}. &\\
\end{array}\\
$$
\end{enumerate}
\end{corollary}
La base $\{X_{0},X_{1},X_{2},\dots,X_{n-1}\}$ ainsi d\'{e}finie est appel\'{e}e base adapt\'{e}e.

\section{Classification des alg\`{e}bres de Lie quasi-filiformes de rang non nul}
Rappelons qu'un tore maximal de d\'{e}rivations de $\frak{g}$ est une sous alg\`{e}bre ab\'{e}lienne maximale de l'alg\`{e}bre des d\'{e}rivations $Der(\frak{g})$ form\'ee de d\'{e}rivations ad-diagonalisables \cite{Ma}. Tous ces tores maximaux sont conjugu\'{e}s par automorphismes, ce qui permet de d\'{e}finir le rang d'une alg\`{e}bre de Lie comme la dimension commune des tores maximaux de $\frak{g}$.
\begin{theorem}
\label{theo:2}
Soit $\frak{g}$ une alg\`{e}bre de Lie quasi-filiforme de dimension $n$ admettant une d\'{e}rivation diagonale $f$ non nulle. Il existe alors une base de $\frak{g}$, $\{Y_{0},Y_{1},\dots,Y_{n-1}\}$, form\'{e}e de vecteurs propres de $f$ dont les crochets v\'{e}rifient l'un des cas suivants:
\begin{enumerate}

\item Si ${\rm gr}\frak{g}\simeq \frak{L_{n-1}}\oplus \mathbb{C} \quad(n\ge 4)$,
\begin{enumerate}
\item $\frak{g}=L_{n-1}\oplus \mathbb{C}$
$$
\lbrack Y_{0},Y_{i}\rbrack=Y_{i+1}, \quad 1 \le i \le n-3,
$$
$$
f\sim diag(\lambda_{0},\lambda_{1},\lambda_{0}+\lambda_{1},2\lambda_{0}+\lambda_{1},\dots,(n-3)\lambda_{0}+\lambda_{1},\lambda_{n-1})
$$
\item \label{1.(b)} $\frak{g}=A_{n-1}^{k}(\alpha_{1},\dots,\alpha_{t-1})\oplus \mathbb{C}, \quad t=\frac{n-k}{2}, \quad 2 \le k \le n-4$
$$
\begin{array}{lll}
\lbrack Y_{0},Y_{i}\rbrack=Y_{i+1}, & 1 \le i \le n-3, &\\
\lbrack Y_{i},Y_{i+1}\rbrack=\alpha_{i} Y_{2i+k}, & 1 \le i \le t-1, &\\
\lbrack Y_{i},Y_{j}\rbrack=a_{i,j} Y_{i+j+k-1}, & 1 \le i <j, & i+j \le n-k-1,\\
\end{array}
$$
$$
f\sim diag(\lambda_{0},k\lambda_{0},(k+1)\lambda_{0},(k+2)\lambda_{0},\dots,(k+n-3)\lambda_{0},\lambda_{n-1})
$$
\item \label{1.(c)} $\frak{g}=L_{n-1} \overrightarrow{\oplus}_{l} \mathbb{C} \quad (2 \le l \le n-3)$
$$
\begin{array}{lll}
\lbrack Y_{0},Y_{i}\rbrack & =Y_{i+1}, & 1 \le i \le n-3,\\
\lbrack Y_{i},Y_{n-1}\rbrack & =Y_{i+l-2}, & 1 \le i \le n-l,\\
\end{array}
$$
$$
f\sim diag(\lambda_{0},\lambda_{1},\lambda_{0}+\lambda_{1},2\lambda_{0}+\lambda_{1},\dots,(n-3)\lambda_{0}+\lambda_{1},l\lambda_{0})
$$
\item \label{1.(d)} $\frak{g}=A_{n-1}^{k}(\alpha_{1},\dots,\alpha_{t-1})\overrightarrow{\oplus}_{l} \mathbb{C} \quad t=\frac{n-k}{2} \quad 2 \le k \le n-4 \quad 2 \le l \le n-3$
$$
\begin{array}{lll}
\lbrack Y_{0},Y_{i}\rbrack=Y_{i+1}, & 1 \le i \le n-3, &\\
\lbrack Y_{i},Y_{n-1}\rbrack=Y_{i+l-2}, & 1 \le i \le n-l, &\\
\lbrack Y_{i},Y_{i+1}\rbrack=\alpha_{i} Y_{2i+k}, & 1 \le i \le t-1, &\\
\lbrack Y_{i},Y_{j}\rbrack=a_{i,j} Y_{i+j+k-1}, & 1 \le i <j, & i+j \le n-k-1,\\
\end{array}
$$
$$
f\sim diag(\lambda_{0},k\lambda_{0},(k+1)\lambda_{0},(k+2)\lambda_{0},\dots,(k+n-3)\lambda_{0},l\lambda_{0})
$$
\end{enumerate}

\item Si ${\rm gr}\frak{g}\simeq \frak{Q_{n-1}}\oplus \mathbb{C} \quad (n\ge 7, \: n \; {\rm impair})$,
\begin{enumerate}
\item $\frak{g}=Q_{n-1}\oplus \mathbb{C}$
$$
\begin{array}{lll}
\lbrack Y_{0},Y_{i}\rbrack=Y_{i+1}, & 1 \le i \le n-4,\\
\lbrack Y_{i},Y_{n-i-1}\rbrack=(-1)^{i-1}Y_{n-2}, & 1\leq i\leq \frac{n-3}{2},\\
\end{array}
$$
$$
f\sim diag(\lambda_{0},\lambda_{1},\lambda_{0}+\lambda_{1},2\lambda_{0}+\lambda_{1},\dots,(n-4)\lambda_{0}+\lambda_{1},(n-4)\lambda_{0}+2\lambda_{1},\lambda_{n-1})
$$
\item \label{2.(b)} $\frak{g}=B_{n-1}^{k}(\alpha_{1},\dots,\alpha_{t-1})\oplus \mathbb{C} \quad t=\frac{n-k-1}{2} \quad 2 \le k \le n-5$
$$
\begin{array}{lll}
\lbrack Y_{0},Y_{i}\rbrack=Y_{i+1}, & 1 \le i \le n-4, &\\
\lbrack Y_{i},Y_{n-i-2}\rbrack =(-1)^{i-1}Y_{n-2}, & 1\leq i\leq \frac{n-3}{2}, &\\
\lbrack Y_{i},Y_{i+1}\rbrack=\alpha_{i} Y_{2i+k-1}, & 1 \le i \le t-1, &\\
\lbrack Y_{i},Y_{j}\rbrack=a_{i,j} Y_{i+j+k-1}, & 1 \le i <j, & i+j \le n-k-2,\\
\end{array}
$$
$$
f\sim diag(\lambda_{0},k\lambda_{0},(k+1)\lambda_{0},(k+2)\lambda_{0},\dots,(n-4+k)\lambda_{0},(n-4+2k)\lambda_{0},\lambda_{n-1})
$$
\item \label{2.(c)} $\frak{g}=Q_{n-1}\overrightarrow{\oplus}_{l}^{a} \mathbb{C}\quad 2 \le l \le n-4$
$$
\begin{array}{ll}
\lbrack Y_{0},Y_{i}\rbrack=Y_{i+1}, & 1 \le i \le n-4,\\
\lbrack Y_{i},Y_{n-i-2}\rbrack=(-1)^{i-1}Y_{n-2}, & 1\leq i\leq \frac{n-3}{2},\\
\lbrack Y_{i},Y_{n-1}\rbrack=Y_{i+l}, & 1\leq i\leq n-l-3,\\
\end{array}
$$
$$
f\sim diag(\lambda_{0},\lambda_{1},\lambda_{0}+\lambda_{1},2\lambda_{0}+\lambda_{1},\dots,(n-4)\lambda_{0}+\lambda_{1},(n-4)\lambda_{0}+2\lambda_{1},l\lambda_{0})
$$
\item \label{2.(d)} $\frak{g}=B_{n-1}^{k}(\alpha_{1},\dots,\alpha_{t-1})\overrightarrow{\oplus}_{l}^{a} \mathbb{C} \quad t=\frac{n-k-1}{2} \quad 2 \le k \le n-5, \; 2 \le l \le n-4$
$$
\begin{array}{lll}
\lbrack Y_{0},Y_{i}\rbrack=Y_{i+1}, & 1 \le i \le n-4, &\\
\lbrack Y_{i},Y_{n-i-2}\rbrack=(-1)^{i-1}Y_{n-2}, & 1\leq i\leq \frac{n-3}{2}, &\\
\lbrack Y_{i},Y_{i+1}\rbrack=\alpha_{i} Y_{2i+k}, & 1 \le i \le t-1, &\\
\lbrack Y_{i},Y_{j}\rbrack=a_{i,j} Y_{i+j+k-1}, & 1 \le i <j, & i+j \le n-k-2,\\
\lbrack Y_{i},Y_{n-1}\rbrack=Y_{i+l}, & 1\leq i\leq n-l-3, &\\
\end{array}
$$
$$
f\sim diag(\lambda_{0},k\lambda_{0},(k+1)\lambda_{0},(k+2)\lambda_{0},\dots,(n-4+k)\lambda_{0},(n-4+2k)\lambda_{0},l\lambda_{0})
$$
\item \label{2.(e)} $\frak{g}=Q_{n-1}\overrightarrow{\oplus}_{l}^{b} \mathbb{C}\quad 2 \le l \le n-4$
$$
\begin{array}{ll}
\lbrack Y_{0},Y_{i}\rbrack=Y_{i+1}, & 1 \le i \le n-4,\\
\lbrack Y_{i},Y_{n-i-2}\rbrack=(-1)^{i-1}Y_{n-2}, & 1\leq i\leq \frac{n-3}{2},\\
\lbrack Y_{0},Y_{n-1}\rbrack=Y_{n-2}, &\\
\lbrack Y_{i},Y_{n-1}\rbrack=Y_{i+l}, & 1\leq i\leq n-l-3,\\
\end{array}
$$
$$
f\sim diag(\lambda_{0},\beta \lambda_{0},(k+1)\lambda_{0},(\beta+2)\lambda_{0},\dots,(n-4+\beta)\lambda_{0},(n-4+2\beta)\lambda_{0},(n-5+2\beta)\lambda_{0})
$$
o\`{u} $\beta=\frac{l-n+5}{2}$
\item \label{2.(f)} $\frak{g}=Q_{n-1}\overrightarrow{\oplus}^{c} \mathbb{C}$
$$
\begin{array}{ll}
\lbrack Y_{0},Y_{i}\rbrack=Y_{i+1}, & 1 \le i \le n-4,\\
\lbrack Y_{i},Y_{n-i-2}\rbrack=(-1)^{i-1}Y_{n-2}, & 1\leq i\leq \frac{n-3}{2},\\
\lbrack Y_{0},Y_{n-1}\rbrack=Y_{n-2}, &\\
\end{array}
$$
$$
f\sim diag(\lambda_{0},\lambda_{1},\lambda_{0}+\lambda_{1},2\lambda_{0}+\lambda_{1},\dots,(n-4)\lambda_{0}\lambda_{1},(n-4)\lambda_{0}+2\lambda_{1},(n-5)\lambda_{0}+2\lambda_{1})
$$
\item \label{2.(g)} $\frak{g}=B_{n-1}^{k}(\alpha_{1},\dots,\alpha_{t-1})\overrightarrow{\oplus}^{c} \mathbb{C} \quad t=\frac{n-k-1}{2} \quad 2 \le k \le n-5$
$$
\begin{array}{lll}
\lbrack Y_{0},Y_{i}\rbrack=Y_{i+1}, & 1 \le i \le n-4, &\\
\lbrack Y_{i},Y_{n-i-2}\rbrack=(-1)^{i-1}Y_{n-2}, & 1\leq i\leq \frac{n-3}{2}, &\\
\lbrack Y_{i},Y_{i+1}\rbrack=\alpha_{i} Y_{2i+k}, & 1 \le i \le t-1, &\\
\lbrack Y_{i},Y_{j}\rbrack=a_{i,j} Y_{i+j+k-1}, & 1 \le i <j, & i+j \le n-k-2,\\
\lbrack Y_{0},Y_{n-1}\rbrack=Y_{n-2}, & &\\
\end{array}
$$
$$
f\sim diag(\lambda_{0},k\lambda_{0},(k+1)\lambda_{0},(k+2)\lambda_{0},\dots,(n-4+k)\lambda_{0},(n-4+2k)\lambda_{0},(n-5+2k)\lambda_{0})
$$
\end{enumerate}

\item Si ${\rm gr}\frak{g}\simeq \frak{L_{n,r}} \quad (n\ge5, \;r\,impair, \;3\le r\le 2[\frac{n-1}{2} ]-1),$
\begin{enumerate}
\item \label{3.(a)} $\frak{g}=\frak{L_{n,r}}$
$$
\begin{array}{ll}
\lbrack Y_{0},Y_{i} \rbrack=Y_{i+1}, &  i=1,\dots,n-3,\\
\lbrack Y_{i},Y_{r-i} \rbrack=(-1)^{i-1}Y_{n-1}, &  i=1,\dots,\frac{r-1}{2}\\
\end{array}
$$
$$
f\sim diag(\lambda_{0},\lambda_{1},\lambda_{0}+\lambda_{1},2\lambda_{0}+\lambda_{1},\dots,(n-3)\lambda_{0}+\lambda_{1},(r-2)\lambda_{0}+2\lambda_{1})
$$
\item \label{3.(b)} $\frak{g}=\frak{C_{n,r}^{k}}(\alpha_{1},\dots,\alpha_{t-1}), \quad 2\le k\le n-4,\;t=[\frac{n-k}{2}]$
$$
\begin{array}{ll}
\lbrack Y_{0},Y_{i} \rbrack=Y_{i+1}, &  i=1,\dots,n-3\\
\lbrack Y_{i},Y_{r-i} \rbrack=\left\{\begin{array}{l}(-1)^{i-1}Y_{n-1}+a_{i,r-i}Y_{r+k-1}\\ (-1)^{i-1}Y_{n-1}\\ \end{array}\right. &  \begin{array}{ll} si\: k\le n-r-1, & i=1,\dots,\frac{r-1}{2}\\ si\: k> n-r-1, & i=1,\dots,\frac{r-1}{2}\\ \end{array}\\
\lbrack Y_{i},Y_{i+1} \rbrack=\alpha_{i}Y_{2i+k}, &  i=1,\dots,t-1\\
\lbrack Y_{i},Y_{j} \rbrack=a_{i,j}Y_{i+j+k-1}, &  1\le i<j<n-1, \;r\ne i+j\le n-k-1,\\
\lbrack Y_{i},Y_{n-1} \rbrack=Y_{2k+r+i-2}, &  i=1,\dots,n-r-2k\\
\end{array}
$$
$$
f\sim diag(\lambda_{0},k\lambda_{0},(1+k)\lambda_{0},(2+k)\lambda_{0},\dots,(n-3+k)\lambda_{0},(r-2+2k)\lambda_{0})
$$
\item \label{3.(c)} $\frak{g}=\frak{D_{n,r}^{k}}, \quad 1\le k\le [\frac{n-r-2}{2}]$
$$
\begin{array}{ll}
\lbrack Y_{0},Y_{i} \rbrack=Y_{i+1}, &  i=1,\dots,n-3\\
\lbrack Y_{i},Y_{r-i} \rbrack=(-1)^{i-1}Y_{n-1}, &  i=1,\dots,\frac{r-1}{2}\\
\lbrack Y_{i},Y_{n-1} \rbrack=Y_{2k+r+i-1}, &  i=1,\dots,n-r-2k-1\\
\end{array}
$$
$$
f\sim diag(\lambda_{0},(k+\frac{1}{2})\lambda_{0},(k+\frac{3}{2})\lambda_{0},(k+\frac{5}{2})\lambda_{0},\dots,(k+\frac{2n-5}{2})\lambda_{0},(r-1+2k)\lambda_{0})
$$
\end{enumerate}

\item Si ${\rm gr}\frak{g}\simeq \frak{Q_{n,r}}\quad (n\ge7,\; n\,impair, \;r\,impair, \;3\le r\le n-4)$,
\begin{enumerate}
\item \label{4.(a)} $\frak{g}=\frak{Q_{n,r}}$
$$
\begin{array}{ll}
\lbrack Y_{0},Y_{i} \rbrack=Y_{i+1}, &  i=1,\dots,n-4\\
\lbrack Y_{i},Y_{r-i} \rbrack=(-1)^{i-1}Y_{n-1}, &  i=1,\dots,\frac{r-1}{2}\\
\lbrack Y_{i},Y_{n-2-i} \rbrack=(-1)^{i-1}Y_{n-2}, &  i=1,\dots,\frac{n-3}{2}\\
\end{array}
$$
$$
f\sim diag(\lambda_{0},\lambda_{1},\lambda_{0}+\lambda_{1},2\lambda_{0}+\lambda_{1},\dots,(n-4)\lambda_{0}+\lambda_{1},(n-4)\lambda_{0}+2\lambda_{1},(r-2)\lambda_{0}+2\lambda_{1})
$$
\item \label{4.(b)} $\frak{g}=\frak{E_{n,r}^{k}}(\alpha_{1},\dots,\alpha_{t-1}), \quad2\le k\le n-5,\;t=[\frac{n-k-1}{2}]$
$$
\begin{array}{ll}
\lbrack Y_{0},Y_{i} \rbrack=Y_{i+1}, &  i=1,\dots,n-4\\
\lbrack Y_{i},Y_{r-i} \rbrack=\left\{\begin{array}{l}(-1)^{i-1}Y_{n-1}+a_{i,r-i}Y_{r+k-1}\\ (-1)^{i-1}Y_{n-1}\\ \end{array}\right. &   \begin{array}{ll}si\, k\le n-r-2, & i=1,\dots,\frac{r-1}{2}\\ si\, k> n-r-2, & i=1,\dots,\frac{r-1}{2}\\ \end{array} \\
\lbrack Y_{i},Y_{n-2-i} \rbrack=(-1)^{i-1}Y_{n-2}, &  i=1,\dots,\frac{n-3}{2}\\
\lbrack Y_{i},Y_{i+1} \rbrack=\alpha_{i}Y_{2i+k}, &  i=1,\dots,t-1\\
\lbrack Y_{i},Y_{j} \rbrack=a_{i,j}Y_{i+j+k-1}, &  1\le i<j<n-1, \;r\ne i+j\le n-k-2\\
\lbrack Y_{i},Y_{n-1} \rbrack=Y_{2k+r+i-2}, &  i=1,\dots,n-r-2k-1\\
\end{array}
$$
$$
f\sim diag(\lambda_{0},k\lambda_{0},(1+k)\lambda_{0},(2+k)\lambda_{0},\dots,(n-4+k)\lambda_{0},(n-4+2k)\lambda_{0},(r-2+2k)\lambda_{0})
$$
\item \label{4.(c)} $\frak{g}=\frak{F_{n,r}^{k}}, \quad 1\le k\le [\frac{n-r-4}{2}]$
$$
\begin{array}{ll}
\lbrack Y_{0},Y_{i} \rbrack=Y_{i+1}, &  i=1,\dots,n-4\\
\lbrack Y_{i},Y_{r-i} \rbrack=(-1)^{i-1}Y_{n-1}, &  i=1,\dots,\frac{r-1}{2}\\
\lbrack Y_{i},Y_{n-2-i} \rbrack=(-1)^{i-1}Y_{n-2}, &  i=1,\dots,\frac{n-3}{2}\\
\lbrack Y_{i},Y_{n-1} \rbrack=Y_{2k+r+i-1}, &  i=1,\dots,n-r-2k-2\\
\end{array}
$$
$$
f\sim diag(\lambda_{0},(k+\frac{1}{2})\lambda_{0},(k+\frac{3}{2})\lambda_{0},(k+\frac{5}{2})\lambda_{0},\dots,(k+\frac{2n-7}{2})\lambda_{0},(n+2k-3)\lambda_{0},(r+2k-1)\lambda_{0})
$$
\end{enumerate}

\item Si ${\rm gr}\frak{g}\simeq \frak{T_{n,n-4}}\quad (n\ge7,\; n\,impair)$,
\begin{enumerate}
\item \label{5.(a)} $\frak{g}=\frak{T_{n,n-4}}$
$$
\begin{array}{ll}
\lbrack Y_{0},Y_{i} \rbrack=Y_{i+1}, &  i=1,\dots,n-5\\
\lbrack Y_{0},Y_{n-3} \rbrack=Y_{n-2}, &\\
\lbrack Y_{0},Y_{n-1} \rbrack=Y_{n-3}, &\\
\lbrack Y_{i},Y_{n-4-i} \rbrack=(-1)^{i-1}Y_{n-1}, &  i=1,\dots,\frac{n-5}{2}\\
\lbrack Y_{i},Y_{n-3-i} \rbrack=(-1)^{i-1} \frac{n-3-2i}{2} Y_{n-3}, &  i=1,\dots,\frac{n-5}{2}\\
\lbrack Y_{i},Y_{n-2-i} \rbrack=(-1)^{i} (i-1) \frac{n-3-i}{2} Y_{n-2},\quad &  i=2,\dots,\frac{n-3}{2}\\
\end{array}
$$
$$
f\sim diag(\lambda_{0},\lambda_{1},\lambda_{0}+\lambda_{1},2\lambda_{0}+\lambda_{1},\dots,(n-5)\lambda_{0}+\lambda_{1},(n-5)\lambda_{0}+2\lambda_{1},(n-4)\lambda_{0}+2\lambda_{1},(n-6)\lambda_{0}+2\lambda_{1})
$$
\item \label{5.(b)} $\frak{g}=\frak{G_{n,r}^{k}}(\alpha_{1},\dots,\alpha_{t-1}),\quad 2\le k\le n-6,\;t=[\frac{n-k-2}{2}]$
$$
\begin{array}{ll}
\lbrack Y_{0},Y_{i} \rbrack=Y_{i+1}, &  i=1,\dots,n-5\\
\lbrack Y_{0},Y_{n-3} \rbrack=Y_{n-2},\\
\lbrack Y_{0},Y_{n-1} \rbrack=Y_{n-3},\\
\lbrack Y_{1},Y_{n-1} \rbrack=Y_{n-2} &si \; k=2\\
\lbrack Y_{i},Y_{n-4-i} \rbrack=(-1)^{i-1}Y_{n-1}, &  i=1,\dots,\frac{n-5}{2}\\
\lbrack Y_{i},Y_{n-3-i} \rbrack=(-1)^{i-1}\frac{n-3-2i}{2}Y_{n-3}, &  i=1,\dots,\frac{n-5}{2}\\
\lbrack Y_{i},Y_{n-2-i} \rbrack=(-1)^{i} (i-1) \frac{n-2-i}{2} Y_{n-2}, &  i=1,\dots,\frac{n-3}{2},\\
\lbrack Y_{i},Y_{i+1} \rbrack=\alpha_{i}Y_{2i+k}, &  i=1,\dots,t-1\\
\lbrack Y_{i},Y_{j} \rbrack=a_{i,j}Y_{i+j+k-1}, &  1\le i<j<n-2,\, i+j\le n-k-3\\
\end{array}
$$
$$
f\sim diag(\lambda_{0},k\lambda_{0},(1+k)\lambda_{0},(2+k)\lambda_{0},\dots,(n-5+k)\lambda_{0},(n-5+2k)\lambda_{0},(n-4+2k)\lambda_{0},(n-6+2k)\lambda_{0})
$$
\end{enumerate}

\item Si ${\rm gr}\frak{g}\simeq \frak{T_{n,n-3}}\quad (n\ge6,\; n\,pair)$,
\begin{enumerate}
\item \label{6.(a)} $\frak{g}=\frak{T_{n,n-3}}$
$$
\begin{array}{ll}
\lbrack Y_{0},Y_{i} \rbrack=Y_{i+1}, &  i=1,\dots,n-4\\
\lbrack Y_{0},Y_{n-1} \rbrack=Y_{n-2},\\
\lbrack Y_{i},Y_{n-3-i} \rbrack=(-1)^{i-1}Y_{n-1}, &  i=1,\dots,\frac{n-4}{2}\\
\lbrack Y_{i},Y_{n-2-i} \rbrack=(-1)^{i-1} \frac{n-2-2i}{2} Y_{n-2}, &  i=1,\dots,\frac{n-4}{2}\\
\end{array}
$$
$$
f\sim diag(\lambda_{0},\lambda_{1},\lambda_{0}+\lambda_{1},2\lambda_{0}+\lambda_{1},\dots,(n-4)\lambda_{0}+\lambda_{1},(n-4)\lambda_{0}+2\lambda_{1},(n-5)\lambda_{0}+2\lambda_{1})
$$
\item \label{6.(b)} $\frak{g}=\frak{H_{n,r}^{k}}(\alpha_{1},\dots,\alpha_{t-1}), \quad2\le k\le n-5,\;t=[\frac{n-k-1}{2}]$
$$
\begin{array}{ll}
\lbrack Y_{0},Y_{i} \rbrack=Y_{i+1}, &  i=1,\dots,n-4\\
\lbrack Y_{0},Y_{n-1} \rbrack=Y_{n-2},\\
\lbrack Y_{i},Y_{n-3-i} \rbrack=(-1)^{i-1}Y_{n-1}, &  i=1,\dots,\frac{n-4}{2}\\
\lbrack Y_{i},Y_{n-2-i} \rbrack=(-1)^{i-1} \frac{n-2-2i}{2} Y_{n-2}, &  i=1,\dots,\frac{n-4}{2}\\
\lbrack Y_{i},Y_{i+1} \rbrack=\alpha_{i}Y_{2i+k}, &  i=1,\dots,t-1\\
\lbrack Y_{i},Y_{j} \rbrack=a_{i,j}Y_{i+j+k-1}, &  1\le i<j<n-2,\, i+j\le n-k-2\\
\end{array}
$$
$$
f\sim diag(\lambda_{0},k\lambda_{0},(1+k)\lambda_{0},(2+k)\lambda_{0},\dots,(n-4+k)\lambda_{0},(n-4+2k)\lambda_{0},(n-5+2k)\lambda_{0})
$$
\end{enumerate}

\item Si ${\rm gr}\frak{g}\simeq \frak{E_{9,5}^{1}}$ alors $\frak{g}\simeq \frak{E_{9,5}^{1}}$
$$
\begin{array}{ll}
\lbrack Y_{0},Y_{i}\rbrack =Y_{i+1}, & i=1,2,3,4,5,6 \\
\lbrack Y_{0},Y_{8}\rbrack =Y_{6}, & \lbrack Y_{2},Y_{8}\rbrack =-3Y_{7},\\
\lbrack Y_{1},Y_{4}\rbrack =Y_{8}, & \lbrack Y_{1},Y_{5}\rbrack =2Y_{6},\\
\lbrack Y_{1},Y_{6}\rbrack =3Y_{7}, & \lbrack Y_{2},Y_{3}\rbrack =-Y_{8},\\
\lbrack Y_{2},Y_{4}\rbrack =-Y_{6}, & \lbrack Y_{2},Y_{5}\rbrack =-Y_{7},\\
\end{array}
$$
$$
f\sim diag(\lambda_{0},\lambda_{0},2\lambda_{0},3\lambda_{0},4\lambda_{0},5\lambda_{0},6\lambda_{0},7\lambda_{0},5\lambda_{0})
$$

\item Si ${\rm gr}\frak{g}\simeq \frak{E_{9,5}^{2}}$ alors $\frak{g}\simeq \frak{E_{9,5}^{2}}$
$$
\begin{array}{ll}
\lbrack Y_{0},Y_{i}\rbrack =Y_{i+1}, &i=1,2,3,4,5,6\\
\lbrack Y_{0},Y_{8}\rbrack =Y_{6}, & \lbrack Y_{2},Y_{8}\rbrack =-Y_{7},\\
\lbrack Y_{1},Y_{4}\rbrack =Y_{8}, & \lbrack Y_{1},Y_{5}\rbrack =2Y_{6},\\
\lbrack Y_{1},Y_{6}\rbrack =Y_{7}, & \lbrack Y_{2},Y_{3}\rbrack =-Y_{8},\\
\lbrack Y_{2},Y_{4}\rbrack =-Y_{6}, & \lbrack Y_{2},Y_{5}\rbrack =Y_{7},\\
\lbrack Y_{3},Y_{4}\rbrack =-2Y_{7}, &\\
\end{array}
$$
$$
f\sim diag(\lambda_{0},\lambda_{0},2\lambda_{0},3\lambda_{0},4\lambda_{0},5\lambda_{0},6\lambda_{0},7\lambda_{0},5\lambda_{0})
$$

\item Si ${\rm gr}\frak{g}\simeq \frak{E_{9,5}^{3}}$ alors $\frak{g}\simeq \frak{E_{9,5}^{3}}$
$$
\begin{array}{ll}
\lbrack Y_{0},Y_{i}\rbrack =Y_{i+1}, &i=1,2,3,4,5,6\\
\lbrack Y_{0},Y_{8} =Y_{6}, & \lbrack Y_{1},Y_{4}\rbrack =Y_{8},\\
\lbrack Y_{1},Y_{5}\rbrack =2Y_{6}, & \lbrack Y_{2},Y_{3}\rbrack =-Y_{8},\\
\lbrack Y_{2},Y_{4}\rbrack =-Y_{6}, & \lbrack Y_{2},Y_{5}\rbrack =2Y_{7},\\
\lbrack Y_{3},Y_{4}\rbrack =-3Y_{7}, &\\
\end{array}\\
$$
$$
f\sim diag(\lambda_{0},\lambda_{0},2\lambda_{0},3\lambda_{0},4\lambda_{0},5\lambda_{0},6\lambda_{0},7\lambda_{0},5\lambda_{0})
$$

\item Si ${\rm gr}\frak{g}\simeq \frak{E_{7,3}}$ alors $\frak{g}\simeq \frak{E_{7,3}}$
$$
\begin{array}{ll}
\lbrack Y_{0},Y_{i}\rbrack =Y_{i+1}, &i=1,2,3,4\\
\lbrack Y_{0},Y_{6}\rbrack =Y_{4}, & \lbrack Y_{2},Y_{6}\rbrack =-Y_{5},\\
\lbrack Y_{1},Y_{2}\rbrack =Y_{6}, & \lbrack Y_{1},Y_{3}\rbrack =Y_{4},\\
\lbrack Y_{1},Y_{4}\rbrack =Y_{5}, &\\
\end{array}\\
$$
$$
f\sim diag(\lambda_{0},\lambda_{0},2\lambda_{0},3\lambda_{0},4\lambda_{0},5\lambda_{0},3\lambda_{0})
$$

\end{enumerate}
Les param\`{e}tres $(\alpha_{1},\dots,\alpha_{t-1})$ v\'{e}rifient les relations polynomiales d\'{e}coulant des identit\'{e}s de Jacobi et les contantes $a_{i,j}$ verifient le syst\`{e}me:
$$
\begin{array}{l}
a_{i,i}=0,\\
a_{i,i+1}=\alpha_{i},\\
a_{i,j}=a_{i+1,j}+a_{i,j+1}.
\end{array}
$$

\end{theorem}

Dans la d\'{e}monstration du th\'{e}or\`{e}me, on souligne les cas o\`{u} ${\rm gr}\frak{g}$ est isomorphe \`{a} $L_{n-1} \overrightarrow{\oplus} \mathbb{C}$ et \`{a} ${\rm gr}\frak{g}\simeq \frak{L_{n,r}}$, pour les autres cas la m\'{e}thode est analogue.\\
\begin{proof}
Supposons que $\frak{g}$ est du type $t_{1}$, alors son alg\`{e}bre gradu\'{e}e est isomorphe \`{a} $\frak{L_{n-1}}\oplus \mathbb{C}$ ou bien \`{a} $\frak{Q_{n-1}}\oplus \mathbb{C}$.\\
Soit $f$ une d\'{e}rivation diagonale de $\frak{g}$ et consid\'{e}rons $\{X_{0},X_{1},\dots,X_{n-1}\}$ une base adapt\'{e}e de $\frak{g}$. On peut trouver trois vecteurs propres de $f$, not\'{e}s $Y_{0},Y_{1},Y_{n-1}$, n'appartenant pas \`{a} l'alg\`{e}bre d\'{e}riv\'{e}e et tels que:
\begin{align*}
Y_{n-1} & = c_{n-1}X_{n-1}+c_{0}X_{0}+c_{1}X_{1}+\sum_{i=2}^{n-2}c_{i}X_{i},\\
Y_{0} & = a_{n-1}X_{n-1}+a_{0}X_{0}+a_{1}X_{1}+\sum_{i=2}^{n-2}a_{i}X_{i},\\
Y_{1} & = b_{n-1}X_{n-1}+b_{0}X_{0}+b_{1}X_{1}+\sum_{i=2}^{n-2}b_{i}X_{i},\\
\end{align*}
avec
$$
{\rm det}\left(
\begin{array}{ccc}
c_{n-1} & a_{n-1} & b_{n-1}\\
c_{0} & a_{0} & b_{0}\\
c_{1} & a_{1} & b_{1}\\
\end{array}
\right)
\ne 0.
$$
On peut alors choisir $Y_{0}$ et $Y_{1}$ de fa\c{c}on que $\Delta={\rm det}\left(
\begin{array}{cc}
a_{0} & b_{0}\\
a_{1} & b_{1}\\
\end{array}
\right)
$ soit non nul.\\
Quand ${\rm gr}\frak{g}$ est isomorphe \`{a} $\frak{L_{n-1}}\oplus \mathbb{C}$, on pose:
$$
Y_{i+1}=\lbrack Y_{0},Y_{i}\rbrack =(a_{0}b_{1}-a_{1}b_{0})X_{i+1}+\sum_{k=i+2}^{n-2}\gamma_{i+1}^{k}X_{k},\quad i=1,\dots,n-3
$$
On en d\'{e}duit que $\{Y_{0},Y_{1},\dots,Y_{n-1}\}$ est une base de vecteurs propres de $f$.
Si $\lambda_{n-1},\lambda_{0},\lambda_{1}$ sont les valeurs propres de $f$ associ\'{e}es respectivement aux vecteurs propres $Y_{n-1},Y_{0},Y_{1}$, alors $Y_{i}$ est un vecteur propre de $f$ avec valeur propre  $(i-1)\lambda_{0}+\lambda_{1}$ pour $i=1,\dots,n-2$.\\
La matrice de changement de variable de la base adapt\'{e}e $\{X_{n-1},X_{0},X_{1},\dots,X_{n-2}\}$ \`{a} la base de vecteurs propres de $f$, $\{Y_{n-1},Y_{0},Y_{1},\dots,Y_{n-2}\}$, est la suivante:
$$
\left(
   \begin{array}{ccclccc}
   c_{n-1} & a_{n-1} & b_{n-1} & \vline &                &        &\\
   c_{0} & a_{0} & b_{0} & \vline &                      & 0      &\\
   c_{1} & a_{1} & b_{1} & \vline &                      &        &\\ \cline{1-3}
   \vdots & \vdots & \vdots & & \Delta&        &\\
         &       &       & & \vdots & \ddots &\\
         &       &       & &        &        & \Delta\\
   \end{array}
\right)
$$
On remarque que le d\'{e}terminant de cette matrice est non nul et on calcule les crochets restants dans la base de vecteurs propres $\{Y_{n-1},Y_{0},Y_{1},\dots,Y_{n-2}\}$.\\
Comme $\lbrack Y_{0},Y_{n-1}\rbrack \in \frak{g}_{3}$ et $f$ est une d\'{e}rivation, on peut supposer que $\lbrack Y_{0},Y_{n-1}\rbrack=0$ en faisant, si n\'{e}cessaire, un changement de variable sur $Y_{n-1}$ .
Par ailleurs, $\lbrack Y_{1},Y_{n-1}\rbrack$ appartient aussi \`{a} $\frak{g}_{3}$ alors:
$$
\lbrack Y_{1},Y_{n-1}\rbrack=\sum_{i=3}^{n-2}d^{i}Y_{i}
$$
En imposant que $f$ est une d\'{e}rivation, on obtient que:
$$
\lbrack (i-1)\lambda_{0}-\lambda_{n-1}\rbrack d^{i}=0 \quad i=3,\dots,n-2
$$
Si on consid\`{e}re en plus les conditions de Jacobi on trouve les cas suivants:
\begin{description}
\item{a)} Lorsque $\lambda_{n-1}=l\lambda_{0}$ pour un certain $l \in \{2,\dots,n-3\}$, on obtient:
$$
\begin{array}{ll}
\lbrack Y_{j},Y_{n-1}\rbrack=Y_{l+j-2}, & 1 \le j \le n-l,\\
\lbrack Y_{j},Y_{n-1}\rbrack=0, & n-l+1 \le j \le n-2.
\end{array}
$$
\item{b)} Si la condition ant\'{e}rieure ne se v\'{e}rifie pas on trouve que:
$$
\lbrack Y_{j},Y_{n-1}\rbrack=0, \quad 1 \le j \le n-2.
$$
\end{description}
Dans ce dernier cas $\frak{g}$ est de la forme $\frak{g^{\prime}} \oplus \mathbb{C}$ o\`{u} $gr(\frak{g^{\prime}})$ est isomorphe \`{a} $L_{n-1}$ puis, en utilisant des arguments similaires, on d\'{e}duit que $\frak{g^{\prime}}$ est isomorphe \`{a} $L_{n-1}$ ou bien elle est du type $A_{n-1}^{k}$.
Par contre, dans le premier cas on obtient que $\frak{g}$ est d'une des formes $L_{n-1} \overrightarrow{\oplus}_{l} \mathbb{C}$ ou $A_{n-1}^{k}(\alpha_{1},\dots,\alpha_{t-1})\overrightarrow{\oplus}_{l} \mathbb{C}$, pr\'{e}cis\'{e}es dans l'\'{e}nonc\'{e} du th\'{e}or\`{e}me.

Quand ${\rm gr}\frak{g}$ est isomorphe \`{a} $\frak{Q_{n-1}}\oplus \mathbb{C}$, le m\^{e}me changement de variable et une d\'{e}monstration analogue nous permettent de d\'{e}duire les r\'{e}sultats.

Supposons  maintenant que $\frak{g}$ est du type $t_{r}$ avec $r\in \{3,\dots,n-2\}$. Si $f$ est une d\'{e}rivation diagonale de $\frak{g}$ et $\{X_{0},X_{1},\dots,X_{n-1}\}$ une base adapt\'{e}e alors on peut trouver deux vecteurs propres de $f$, d\'{e}not\'{e}s par $Y_{0}$ et $Y_{1}$, qui n'appartiennent pas \`{a} l'alg\`{e}bre d\'{e}riv\'{e}e et v\'{e}rifiant:
\begin{align*}
Y_{0} & = a_{0}X_{0}+a_{1}X_{1}+\sum_{i=2}^{n-1}a_{i}X_{i},\\
Y_{1} & = b_{0}X_{0}+b_{1}X_{1}+\sum_{i=2}^{n-1}b_{i}X_{i},\\
\end{align*}
avec
$$
\Delta={\rm det}\left(
\begin{array}{cc}
a_{0} & b_{0}\\
a_{1} & b_{1}\\
\end{array}
\right)
\ne 0.
$$
On peut supposer que $a_{0}=1$.\\
Quand ${\rm gr}\frak{g}$ est isomorphe \`{a} $\frak{L_{n,r}}$, on pose:
\begin{eqnarray*}
Y_{i+1}&=&\lbrack Y_{0},Y_{i}\rbrack \quad i=1,\dots,n-3\\
Y_{n-1}&=&\lbrack Y_{1},Y_{r-1}\rbrack
\end{eqnarray*}
et donc:
\begin{eqnarray*}
Y_{i}&=&\Delta \, X_{i}+\sum_{k=i+1}^{n-1}\gamma_{i}^{k}X_{k},\quad i=2,\dots,r-1,\\
Y_{i}&=&\Delta \, X_{i}+\sum_{k=i+1}^{n-2}\gamma_{i}^{k}X_{k},\quad i=r+1,\dots,n-2,\\
Y_{r}&=&\Delta (X_{r}+a_{1}X_{n-1})+\sum_{k=r+1}^{n-2}\gamma_{r}^{k}X_{k},\\
Y_{n-1}&=&\Delta (b_{0}X_{r}+b_{1}X_{n-1})+\sum_{k=r+1}^{n-2}\gamma_{r}^{k}X_{k}.\\
\end{eqnarray*}
La matrice de changement de base de $\{X_{0},X_{1},\dots,X_{r-1},X_{r},X_{n-1},X_{r+1},\dots,X_{n-2}\}$ \`{a} $\{Y_{0},Y_{1},\dots,Y_{r-1},Y_{r},Y_{n-1},Y_{r+1},\dots,Y_{n-2}\}$ est la suivante:
$$
\left(
   \begin{array}{ccccccccc}
        1 & b_{0}  &        &        &              &             &        &        &\\
    a_{1} & b_{1}  &        &        &              &             &    0   &        &\\
   \vdots & \vdots & \Delta &        &              &             &        &        &\\
          &        & \vdots & \ddots &              &             &        &        &\\
          &        &        &        & \Delta       & \Delta b_{0}&        &        &\\
          &        &        &        & \Delta a_{1} & \Delta b_{1}&        &        &\\
          &        &        &        & \vdots       & \vdots      & \Delta &        &\\
          &        &        &        &              &             & \vdots & \ddots &\\
          &        &        &        &              &             &        &        & \Delta\\
   \end{array}
\right)
$$
$\{Y_{0},Y_{1},\dots,Y_{n-1}\}$ est une base de vecteurs propres de $f$, si $\lambda_{0}$ et $\lambda_{1}$ repr\'{e}sentent respectivement les valeurs propres associ\'{e}es aux vecteurs $Y_{0}$ et $Y_{1}$ alors, dans cette base, $f$ est de la forme:
$$
f\sim diag(\lambda_{0},\lambda_{1},\lambda_{0}+\lambda_{1},2\lambda_{0}+\lambda_{1},\dots,(n-3)\lambda_{0}+\lambda_{1},(r-2)\lambda_{0}+2\lambda_{1})
$$
On calcule les crochets dans cette nouvelle base:
$$
\lbrack Y_{0}, Y_{n-2} \rbrack =0,
$$
$$
\lbrack Y_{0}, Y_{n-1} \rbrack =\left\{
\begin{array}{ll}
\Delta b_{0}Y_{r+1}+\sum_{k=2}^{n-2}c_{k}Y_{r+k} & si\,r\ne n-2\\
0 & si\,r=n-2\\
\end{array}
\right..
$$
Si $r\ne n-2$, on impose que $f$ est une d\'{e}rivation:
$$
f(\lbrack Y_{0}, Y_{n-1} \rbrack)=\lbrack (r-1)\lambda_{0}+2\lambda_{1}\rbrack \lbrack Y_{0}, Y_{n-1} \rbrack \Rightarrow \left\{
\begin{array}{ll}
b_{0}(\lambda_{1}-\lambda_{0})=0\\
c_{k}(\lambda_{1}-k\lambda_{0})=0 & pour\,k=2,\dots,n-r-2\\
\end{array}
\right.
$$
Si $2\le i \le r-1$:
$$
\lbrack Y_{1}, Y_{i} \rbrack = b_{0}Y_{i+1}+\sum_{k=2}^{n-i-1}d_{k}Y_{i+k}
\Rightarrow \left\{
\begin{array}{ll}
b_{0}(\lambda_{1}-\lambda_{0})=0\\
d_{k}(\lambda_{1}-k\lambda_{0})=0 & pour\,k=2,\dots,n-i-2\\
d_{n-i-1}\lambda_{0}=0\\
\end{array}
\right.
$$
Si $r\le i \le n-3$:
\begin{eqnarray*}
\lbrack Y_{1}, Y_{i} \rbrack \quad & = & b_{0}Y_{i+1}+\sum_{k=2}^{n-i-2}d_{k}Y_{i+k}
\Rightarrow \left\{
\begin{array}{ll}
b_{0}(\lambda_{1}-\lambda_{0})=0\\
d_{k}(\lambda_{1}-k\lambda_{0})=0 & pour\,k=2,\dots,n-i-2\\
\end{array}
\right.\\
\lbrack Y_{1}, Y_{n-2} \rbrack & = & 0\\
\lbrack Y_{1}, Y_{n-1} \rbrack & = & b_{0}^{2}Y_{r+1}+\sum_{l=3}^{n-r-1}e_{l}Y_{l+r-1}
\Rightarrow \left\{
\begin{array}{ll}
b_{0}(\lambda_{1}-\lambda_{0})=0\\
d_{l}(2\lambda_{1}-l\lambda_{0})=0 & pour\,l=3,\dots,n-r-1\\
\end{array}
\right.\\
\end{eqnarray*}
Si $2\le i \le \frac{r-1}{2}$:
$$
\lbrack Y_{i}, Y_{r-i} \rbrack = (-1)^{i-1}(Y_{n-1}-b_{0}Y_{r})+\sum_{k=2}^{n-r-1}g_{k}Y_{r+k-1}
\Rightarrow \left\{
\begin{array}{ll}
b_{0}(\lambda_{1}-\lambda_{0})=0\\
g_{k}(\lambda_{1}-k\lambda_{0})=0 & pour\,k=2,\dots,n-r-1\\
\end{array}
\right.
$$
Si $1<i<j<n-1$ et $i+j\le r-1$:
$$
\lbrack Y_{i}, Y_{j} \rbrack = \sum_{k=2}^{n-i-j}h_{k}Y_{k+i+j-1}
\Rightarrow \left\{
\begin{array}{ll}
h_{k}(\lambda_{1}-k\lambda_{0})=0 & pour\,k=2,\dots,n-i-j-1\\
h_{n-i-j}\lambda_{0}=0\\
\end{array}
\right.
$$
Si $1<i<j<n-1$ et $i+j\ge r+1$:
$$
\lbrack Y_{i}, Y_{j} \rbrack = \sum_{k=2}^{n-i-j-1}h_{k}Y_{k+i+j-1}
\Rightarrow
h_{k}(\lambda_{1}-k\lambda_{0})=0 \, pour\,k=2,\dots,n-i-j-1
$$
Si $1<i<n-1$:
$$
\lbrack Y_{i}, Y_{n-1} \rbrack = \sum_{l=3}^{n-i-r}p_{l}Y_{l+i+r-2}
\Rightarrow
p_{l}(2\lambda_{1}-l\lambda_{0})=0 \, pour\,l=3,\dots,n-i-r
$$

Ces calculs nous font envisager diff\'{e}rentes possibilit\'{e}s, en consid\'{e}rant, en plus, les conditions de Jacobi.
\begin{enumerate}
\item Si $\lambda_{1}=\lambda_{0}$, il suffit de faire le changement de variable $Y_{1}\to Y_{1}-b_{0}Y_{0},Y_{n-1}\to Y_{n-1}-b_{0}Y_{r}$ pour trouver l'alg\`{e}bre $\frak{L_{n,r}}$
\item Si $\lambda_{0}=0$, on peut supposer $\lambda_{1}=1$ et il existe un changement de variable avec lequel on obtient l'alg\`{e}bre $\frak{L_{n,r}}$.
\item Si $\lambda_{1}=k\lambda_{0}$ avec $2\le k\le n-4$, on obtient une alg\`{e}bre du type $\frak{C_{n,r}^{k}}(\alpha_{1},\dots,\alpha_{t-1})$ o\`{u} $t=[\frac{n-k}{2}]$.
\item Si $2\lambda_{1}=(2k+1)\lambda_{0}$ avec $1\le k\le [\frac{n-r-2}{2}]$, on obtient une alg\`{e}bre du type $\frak{g}=\frak{D_{n,r}^{k}}$.
\end{enumerate}
Dans le reste des cas, le m\^{e}me changement de variable nous permet d'obtenir les familles d'alg\`{e}bres d\'{e}crites dans le th\'{e}or\`{e}me ~\ref{theo:2}.
\end{proof}
\begin{remark}
Sauf isomorphisme, il y a seulement deux alg\`{e}bres quasi-filiformes de rang $3$, respectivement isomorphes \`{a} $\frak{L_{n-1}}\oplus \mathbb{C}$ et \`{a} $\frak{Q_{n-1}}\oplus \mathbb{C}$.
\par
Les alg\`{e}bres quasi-filiformes de rang $2$ sont $\frak{L_{n,r}}$, $\frak{Q_{n,r}}$, $\frak{T_{n,n-4}}$, $\frak{T_{n,n-3}}$ ou bien de la forme ~\ref{1.(b)},~\ref{1.(c)},~\ref{2.(b)},~\ref{2.(c)},~\ref{2.(f)},~\ref{2.(g)} d'apr\`{e}s la terminologie du th\'{e}or\`{e}me ~\ref{theo:2}.
\par
Les alg\`{e}bres quasi-filiformes de rang $1$ sont $\frak{E_{9,5}^{1}}$, $\frak{E_{9,5}^{2}}$, $\frak{E_{9,5}^{3}}$, $\frak{E_{7,3}}$ ou bien correspondent \`a un des alg\`ebres du type ~\ref{1.(d)},~\ref{2.(d)},~\ref{2.(e)},~\ref{3.(b)},~\ref{3.(c)},~\ref{4.(b)},~\ref{4.(c)},~\ref{5.(b)},~\ref{6.(b)} d\'ecrites dans le th\'{e}or\`{e}me ~\ref{theo:2}.\\
\end{remark}

\subsection*{Remerciements}
L'auteur est soutenue par le projet de recherche MTM2006-09152 du Ministerio de Educaci\'{o}n y Ciencia, et remercie aussi la Fundaci\'{o}n Ram\'{o}n Areces pour financer sa bourse pr\'{e}doctorale.

\end{document}